\documentclass[12pt]{article}
\usepackage[all]{xy}
\usepackage{amsmath}
\usepackage{amssymb}
\usepackage{amsthm}

\newtheorem{theorem}{Theorem}
\newtheorem{lemma}[theorem]{Lemma}
\newtheorem{proposition}[theorem]{Proposition}
\newtheorem{corollary}[theorem]{Corollary}

\newtheorem{fact}{Fact}
\theoremstyle{definition}
\newtheorem{definition}[theorem]{Definition}
\newtheorem{property}[fact]{Property}

\begin{document}
\newcommand{\xx}{\mathbf{x}}
\newcommand{\yy}{\mathbf{y}}
\newcommand{\conj}{\mathbf{conj}}
\newcommand{\hypfails}{\mathbf{hypfails}}
\newcommand{\concholds}{\mathbf{concholds}}
\newcommand{\inl}{\mathbf{inl}}
\newcommand{\inr}{\mathbf{inr}}
\newcommand{\inlinr}{\mathbf{inlinr}}
\newcommand{\reasons}{\mathrm{Reasons}}
\newcommand{\reasonsinf}{\mathrm{Reasons}^\infty}
\newcommand{\meet}{\wedge}	
\newcommand{\rank}{\mathrm{rank}}
\newcommand{\N}{\mathbb{N}}
\newcommand{\0}{\mathbf{0}}
\renewcommand{\subset}{\subseteq}
\newcommand{\straightforward}{\begin{proof} Straightforward.\end{proof}}
\newcommand{\lev}{\mathrm{Lev}}
\bibliographystyle{amsplain}
\title{The Order-Theoretic Structure of Free Heyting Algebras}
\author{Michael O'Connor\\
Department of Mathematics \\Cornell University \\Ithaca NY 14853\\
\texttt{oconnor@math.cornell.edu}}

\maketitle

\begin{abstract}We find an order-theoretic characterization of the Lindenbaum algebra of
intuitionistic propositional logic in $n$ variables.
\end{abstract}

\section{Introduction}
Intuitionistic logic has been explored for many years as a
language for computer science, with a guiding principle being the Brouwer-Heyting-Kolmogorov
interpretation, under which intuitionistic proofs of implication 
are functions and existence proofs require witnesses. Higher-order intuitionistic systems
which can express a great deal of mathematics, such as Girard's System
F and Martin-L\"of's type theory (good references are~\cite{girardpat} and~\cite{beeson}), have been developed and implemented by
prominent computer scientists such as Constable, Huet and Coquand (see~\cite{nuprl} and~\cite{coq}). With
all this development and with the existence of well-established topological, Kripke, and
categorical semantics for intuitionistic systems, it may come as a
surprise that many fundamental structural properties of intuitionistic
propositional calculus have not been developed.
By way of contrast, corresponding issues for classical logics have been settled for
at least 75 years. 

Heyting algebras are an equationally defined class of algebras with operations 
$\vee$, $\wedge$, and $\rightarrow$ and constants $\bot$ and $\top$ (representing ``or,'' ``and,'' ``implies,'' ``false,'' and ``true'' respectively) that stand in the same relation to intuitionistic propositional
logic that Boolean algebras do to classical propositional logic. What follows is a very brief introduction to free Heyting
algebras and a summary of the results that will be presented in this paper.

\newcommand{\iv}{\vdash_{\mathcal{I}}}
\newcommand{\cv}{\vdash_{\mathcal{C}}}

For each $n\in\N$, let $V_n = \{x_1,\ldots, x_n\}$ and let $F_n$ be the set of
propositional sentences in variables $V_n$. 
Let $\simeq^n_i$ and $\simeq^n_c$ be the intuitionistic and classical logical equivalence
relations respectively.  

The classical Lindenbaum algebra $B_n$ is then defined as $F_n/\simeq^n_c$ 
and the intuitionistic Lindenbaum algebra $H_n$ is defined as $F_n/\simeq^n_i$.
The operations $\meet$ and $\vee$ and the constants $\top$ and $\bot$ are naturally defined on $B_n$ and the operations
$\meet$, $\vee$, and $\rightarrow$ and the constants $\top$ and $\bot$ are naturally defined on $H_n$. $B_n$
is then isomorphic to the free Boolean algebra on $n$ generators and $H_n$ is the free Heyting algebra
on $n$ generators. As usual, the order $\leq$ may be defined from $\wedge$ (or from $\vee$). Like all Heyting algebras, each $H_n$ is also a distributive lattice.

%
%
%

The structure of each $B_n$ and of $B_\omega$ is well understood. However, among the free Heyting algebras,
up to this point only $H_1$ has been completely understood. It is known from~\cite{nishimura} that if we let $\phi_1 = \neg x_1$,
$\psi_1 = x_1$, $\phi_{i+1} = \phi_i \rightarrow \psi_i$, and $\psi_{i+1} = \phi_i \vee \psi_i$, then
each propositional formula in the single variable $x_1$ is intuitionistically equivalent to exactly one formula
in $\{\bot\}\cup  \{\phi_i\mid i \in\omega\} \cup  \{\psi_i\mid i \in\omega\}\cup\{\top\}$.
Further, we can easily write down conditions characterizing the order on those formulas, so that the structure of $H_1$ is completely characterized.

Although not a complete list, the reader is referred to~\cite{bellissima}, \cite{urq}, \cite{ghilardi}, \cite{butz},
and \cite{darniere} for previous work on free or finitely generated Heyting algebras.  A very useful construction is contained in~\cite{bellissima} which will we avail ourselves of in this paper and which is described in Section~\ref{notation} below.

It is known that to characterize $H_n$ 
 it suffices to characterize the join-irreducible elements of $H_n$ 
 (see~\cite{urq}), as every element of $H_n$ 
 is equivalent to a unique join of join-irreducibles. Let $J_n$ be the poset of join-irreducible elements of $H_n$ 
.

In this paper we will characterize each $J_n$. 

\section{The Main Theorem}

\begin{definition}[Quasisemilattices]
A quasisemilattice (qsl) is a poset  $(P,\leq)$ such that for any two elements
$p$ and $q$, the set $\{r\in P \mid r\leq p, r\leq q\}$ of lower bounds of $p$ and $q$ has only finitely many
maximal elements and such that every lower bound of $p$ and $q$ is below at least one such maximal element. 

A qsl is called bounded (a bqsl) if it has a minimum element.

A qsl is called locally finite if it is locally finite under the relation 
$R(p,q,r)$ which holds iff $r$ is a maximal lower bound of $p$ and $q$.

A bqsl embedding between two bqsls $Q_1$ and $Q_2$
is an order-embedding that respects the minimal element $\bot$ and respects, for each pair of elements
$p,q\in Q_1$ the number and identity of the maximal elements of $\{r\in Q_1\mid
r\leq p,r\leq q\}$.

A qsl is universal countable bounded locally finite if it is countable, bounded, and locally
finite, and embeds all countable, bounded, locally finite qsls.
\end{definition}
By a standard Fra\"iss\'e argument, there is a unique universal countable homogeneous locally finite bounded quasisemilattice.  Let it be $Q$. Let $Q'$ be $Q$ with a maximum element added and
the minimum element removed.

\begin{definition}[$J_{1,n}$, $J_{2,n}$, $J_{3,n}$]
For each $n$, let $J_{1,n} = \{\phi\in J_n \mid (\exists^{<\infty}\psi)\,\psi<\phi\}$.

Let $J_{2,n}$ be the set of minimal elements of $J_n - J_{1,n}$.

Let $J_{3,n} = J_n - (J_{1,n}\cup J_{2,n})$.
\end{definition}
As discussed below in Section~\ref{j1n}, $J_{1,n}$ is characterized completely in~\cite{bellissima}.
\begin{theorem} For all $n\geq 2$, $J_{2,n}$ is a countably infinite antichain and
$J_{3,n}$ is isomorphic to $Q'$.

Every element of $J_{2,n}$ has an element of $J_{1,n}$ below it. Every element of 
$J_{3,n}$ has an element of $J_{2,n}$ below it.

If $x\in J_{i,n}$ and $y\in J_{j,n}$ and $x\leq y$, then $i\leq j$.
\end{theorem}
\begin{proof}
Proposition~\ref{pj1n} states that every element of $J_n$ has an element of $J_{1,n}$
below it. 

Proposition~\ref{j2} states that $J_{2,n}$ is a countably infinite antichain and that
every element of $J_{3,n}$ has an element of $J_{2,n}$ below it.

Proposition~\ref{j3} states that $J_{3,n}$ is isomorphic to $Q'$.

The final statement follows immediately from the definitions.
\end{proof}

\newcommand{\node}{\mathop{\text{node}}}

\section{Terminology, Notation, and Bellissima's Construction}\label{notation}
This paper will make heavy use of Bellissima's construction from~\cite{bellissima}.  To make this paper self-contained, the construction and the relevant facts about it will be stated here.

If $\alpha < \beta$ are nodes in some Kripke model, then $\beta$ is called a successor
of $\alpha$ and $\alpha$ a predecessor of $\beta$.  If there is no $\gamma$
with $\alpha \leq \gamma\leq \beta$, then $\beta$ is called an immediate successor of 
$\alpha$. The assertions ``$\alpha$ is above $\beta$'' and ``$\beta$ is below $\alpha$''
both mean $\alpha \geq \beta$.

For any node $\alpha$ in some Kripke model defined over some set $V$ of propositional variables,
we let $w(\alpha) = \{v\in V\mid \alpha\Vdash v\}$. 

For any formula $\phi$ where some Kripke model $K$ is given by context, let $k(\phi) = 
\{\alpha\in K\mid \alpha\Vdash \phi\}$.

In~\cite{bellissima}, for each $n \geq 1$, Bellissima defines (explicitly) a Kripke model $K_n$ with the following properties: 

\begin{property}For every set of nodes $S\subset K_n$ such that $S$ has at least two minimal elements (under the partial order $K_n$) and every set of atomic formulas
$U \subset \bigcap_{\alpha \in S} w(\alpha)$, there is a unique node we will denote $\node(S,U)$ whose immediate successors are exactly the minimal elements of $S$ and such that $w(\node(S,U)) = U$.
\end{property}

\begin{property} For every set of node $\alpha\in K_n$ and every set of atomic formulas $U \subsetneq w(\alpha)$,
there is a unique node we will denote $\node(\{\alpha\},U)$ with $\alpha$ as its single immediate successor
and such that $w(\node(\{\alpha\},U)) = U$.
\end{property}

As a notational matter, if $S$ has at least two minimal elements, 
then $\node(S)$ will mean $\node(S,\bigcap_{\alpha\in S} w(\alpha))$,
and if $S$ has a single minimum element $\alpha$, then $\node(S)$ will
denote $\alpha$.

Note that in either case, for all sets of nodes $S$ and nodes $\alpha$, 
$\alpha \geq \node(S)$ iff $\alpha = \node(S)$ or for some $\beta\in S$, $\alpha\geq\beta$.

\begin{property} There is a partition of $K_n$ into levels $L^1_n, L^2_n,\ldots$. Each $L^m_n$ is finite. For all nodes
$\alpha$,  let
$\lev(\alpha)$ be the unique $m$ such that $\alpha \in L^m_n$.
For $m \geq 2$
and all nodes $\alpha$, if $\lev(\alpha) \geq 2$, then $\alpha$ has successors and
$\lev(\alpha) = \max\{\lev(\beta) \mid \beta \geq\alpha\} + 1$. If $\lev(\alpha) = 1$,
then $\alpha$ has no successors.

A node $\alpha$ is of the form
$\node(S,U)$ iff $\lev(\alpha)\geq 2$.
\end{property}

Note that if $\alpha < \beta$, $\lev(\alpha) > \lev(\beta)$.  Also, note that any node has only finitely  many successors.

\begin{property} For $\phi$, $\psi\in F_n$, $\phi \vdash \psi$ iff $k(\phi)\subset k(\psi)$.
\end{property}
\begin{property} For each node $\alpha\in K_n$, there is a $\phi_\alpha$ such that
$k(\phi_\alpha) = \{\beta\in K_n\mid \beta \geq \alpha\}$ and there is a $\phi'_\alpha$ such that 
$k(\phi'_\alpha) = \{\beta\in K_n\mid \beta\not\leq \alpha\}$.
\end{property}
\begin{property} If $n\geq 2$, $|L^1_n|\geq 3$.
\end{property}
\begin{property}\label{reduce} For any node $\alpha$ in any finite Kripke model $K$ over $V_n$, there is a node $\beta\in K_n$ such
that for all formulas $\phi\in F_n$, $\alpha \Vdash \phi$ iff $\beta \Vdash \phi$.
\end{property}
The following facts will be used below and follow without much difficulty directly from the above properties.
\begin{fact}\label{biglevel} For $n\geq 2$ and $m\geq 0$, $|L_n^{m+1}| > |L_n^{m}|$. In particular, there are arbitrarily large levels of
$K_n$.
\end{fact}
The following fact is a more general version of the preceding fact.
\begin{fact}\label{bignotnot} Let $S\subset K_n$, $|S| \geq 3$ and let each element of $S$ be at the same level. Let $S'$
be the downward closure of $S$. Then $|S'\cap L_n^{m+1}| > |S'\cap L_n^m|$ for any $m$ greater than or equal
to the common level of the elements of $S$. 
\end{fact}

\section{$J_{1,n}$}\label{j1n}
Let $P_{1,n}$ be the underlying partial order of $K_n$ with the ordering reversed. 

\begin{proposition}[Implicit in \cite{bellissima}]\label{pj1n} $P_{1,n}$ and $J_{1,n}$ are order-isomorphic.  Every element of $H_n$ besides the minimal element has an element of
$J_{1,n}$ below it. A join-irreducible formula $\phi$ is in $J_{1,n}$ iff $k(\phi)$ is finite.
\end{proposition}
The isomorphism sends $\alpha$ to $\phi_\alpha$.

\section{Join-Irreducibles}
Here we will collect some useful lemmas and propositions.

\newcommand{\subforms}{\mathrm{Subform}}
\newcommand{\constypes}{\mathrm{ConsistentTypes}}
\newcommand{\type}{\mathop{\mathrm{Type}}}
\begin{definition}For any formula $\phi$, let $\subforms(\phi)$ be the set of all subformulas of $\phi$ and let
\[ C(\phi) = \{ T\subset \subforms(S) \mid \phi\in T\text{ and }(\exists \alpha\in K_n)\, (\forall \psi\in \subforms(\phi))\,
\alpha \Vdash \psi\text{ iff } \psi\in T\}\]

For any $\alpha\in K_n$, let $\type_\phi(\alpha) = \{ \psi\in \subforms(\phi) \mid \alpha \Vdash \psi\}$. We will have $\type_\phi(\alpha) \in C(\phi)$ iff $\alpha\Vdash \phi$.
\end{definition}

Note that by Property~\ref{reduce}, we could equivalently replace ``$\exists \alpha\in K_n$''
in the definition of $C(\phi)$ by ``$\exists K\exists \alpha\in K$;'' i.e., we could
allow $\alpha$ to range over all nodes of all finite Kripke models.

\begin{lemma}\label{hypfails}Let $\phi$ be a formula. 
Suppose $\alpha$ is a node in a Kripke model with immediate successors
$\beta_1,\ldots,\beta_m$. If $\{\type_\phi(\beta_i)\mid 1\leq i \leq m\}$ has a minimum
element $T$ under $\subset$, and if for atomic formulas $x_i$, 
$\alpha \Vdash x_i$ iff $x_i\in T$, then for all formulas $\psi\in \subforms(\phi)$,
$\alpha \Vdash \psi$ iff $\psi\in T$.
\end{lemma}
\begin{proof}
By induction on the structure of the formula $\psi$.  For atomic formulas it is assumed.
The $\wedge$ and $\vee$ cases are trivial.

Assume that $\rho\rightarrow \chi\in T$. Then if $\rho\in T$, then $\chi\in T$, so by induction,
if $\alpha \Vdash \rho$ then $\alpha\Vdash \chi$. Since $T$ is minimal, all nodes strictly above
$\alpha$ must force $\rho\rightarrow\chi$. Thus $\alpha\Vdash \rho\rightarrow\chi$.

If $\rho\rightarrow\chi\notin T$, then $\alpha$ cannot force $\rho\rightarrow\chi$, since
there is a node not forcing $\rho\rightarrow\chi$ above it.
\end{proof}

\begin{proposition}For any formula $\phi$, $\phi$ is join-irreducible
iff $C(\phi)$ has a minimum element under $\subset$.
\end{proposition}
\begin{proof}
$\Rightarrow$: We will prove the contrapositive.
Let $\phi$ be a set of formulas such that $C(\phi)$
has minimal elements $\{T_1,\ldots, T_m\}$ with $m\geq 2$.  For $1\leq i \leq m$, let
$\phi_i = \bigwedge_{\psi\in T_i}\psi$.  

Then $\bigvee_{1\leq i\leq m} \phi_i$ is equivalent to $\phi$: every node $\alpha\in
k(\phi)$ is in $k(\phi_i)$ where $i$ is such that $\type_\phi(\alpha)\supset T_i$. Conversely, $\phi$ is a
conjunct of each $\phi_i$.

But no $\phi_i$ is equivalent to $\phi$: any node $\alpha$ with $\type_\phi(\alpha) = T_j$
for $j\ne i$ is in $k(\phi)$ but cannot be in $k(\phi_i)$ since it will not force one of the conjuncts of $\phi_i$. Thus $\phi$ is not join-irreducible.

$\Leftarrow$: Suppose $C(\phi)$ has  a minimum element $T$
under $\subset$ and suppose that $\phi$ is equivalent to $\psi\vee\rho$ 
and not equivalent to either $\psi$ or $\rho$.  Then there must be some
$\alpha\in k(\phi)$ such that $\alpha\Vdash \psi$ and $\alpha\not\Vdash\rho$
and some $\beta\in k(\phi)$ such that $\beta\Vdash \psi$ and $\beta\not\Vdash \rho$.

Let $\gamma\in k(\phi)$ be such that $\type_\phi(\gamma) = T$.  Then $\node(\{\alpha,
\beta,\gamma\})$ is in $k(\phi)$ by Lemma~\ref{hypfails} and is below $\alpha$ and $\beta$.
Thus it can be in neither $k(\psi)$ nor $k(\rho)$ which is a contradiction.
\end{proof}
There is a corollary to this which is interesting in its own right.
\begin{corollary}\label{subformulas}If $\phi$ is not join-irreducible, then it is equivalent to
$\phi_0\vee\cdots\vee\phi_r$, where each $\phi_i$ is a conjunction of subformulas
of $\phi$, each $\phi_i$ is join-irreducible, and no $\phi_i$ is equivalent to $\phi$.
\end{corollary}
\begin{proof}
Let $\phi$ not be join-irreducible. Then $C(\phi)$ has
more than one minimal element.  Say its minimal elements are $T_1,\ldots, T_r$.
Let $\phi_i$ for $1\leq i \leq r$ be $\bigwedge_{\psi\in T_i}\psi$.

Then $\phi\vdash \phi_1\vee\cdots\vee\phi_r$: let $\alpha\in k(\phi)$ and let 
$T_i\subset \type_\phi(\alpha)$. Then $\alpha\in k(\phi_i)$.

Conversely, $\phi_1\vee\cdots\vee\phi_r\vdash \phi$, since $\phi$ is a conjunct of each
$\phi_i$.

No $\phi_i$ implies $\phi$: pick an $\alpha$ such that $\type_\phi(\alpha)=T_j$
for $j\ne i$. Then $\alpha\in k(\phi)-k(\phi_i)$.

Finally, each $\phi_i$ is join-irreducible: Suppose $\psi\vee \rho$ is equivalent to
$\phi_i$, and $\phi_i$ is equivalent to neither $\psi$ nor $\rho$. Then let 
$\alpha_\psi \in k(\phi_i) - k(\psi)$, $\alpha_\rho \in k(\phi_i) - k(\rho)$,
and let $\beta$ be such that $\type_\phi(\beta) = T_i$.  Then $\type_\phi(\node
(\{\beta,\alpha_\psi,\alpha_\rho\})) = T_i$, but $\node(\{\beta,\alpha_\psi,\alpha_\rho\})$
is in neither $k(\psi)$ nor $k(\rho)$.
\end{proof}
\begin{lemma}\label{generaljoin}Suppose $\phi_1,\ldots,\phi_m$ are any formulas, $\psi$ is join-irreducible, and 
$\psi\leq \phi_1\vee\cdots\vee\phi_m$.  Then there is an $i$ such that $\psi\leq \phi_i$
\end{lemma}
\begin{proof}
Since $\psi\leq \phi_1\vee\cdots\vee\phi_m$, $\psi$ is equivalent to 
$\psi\wedge(\phi_1\vee\cdots\vee\phi_m)$ and thus to 
$(\psi\wedge\phi_1)\vee\cdots \vee(\psi\wedge\phi_m)$. Since $\psi$ is join-irreducible,
it must be equivalent to some $\psi\wedge\phi_i$, and therefore must be less than
$\phi_i$.
\end{proof}
\newcommand{\mintype}{\mathop{\text{mintype}}}
\begin{definition}If $\phi$ is a join-irreducible formula, let $\mintype(\phi)$ be
the minimum element of $C(\phi)$
under $\subset$.
\end{definition}

\section{$J_{2,n}$}

\begin{proposition}\label{j2}$J_{2,n}$ is a countably infinite antichain.  Every element of $J_n - (J_{2,n}\cup J_{1,n})$ has an element of $J_{2,n}$ below it. 

An element $\phi$ of $J_n$ is in
$J_{2,n}$ iff $k(\phi)$ is infinite and for all but finitely many $m$, $|L^m_n\cap k(\phi)| = 2$.
\end{proposition}
\begin{proof}
We begin with a few definitions.
\begin{definition}[Well-Positioned Triplets of Nodes]
Let $\alpha$, $\beta$, and $\gamma$ be three distinct nodes in $K_n$.
The ordered triplet $(\alpha,\beta,\gamma)$ is called well-positioned if the following properties hold:

$\lev(\alpha) + 1 = \lev(\beta) + 1 = \lev(\gamma)$. $\gamma < \beta$.
$\gamma \not < \alpha$. $w(\alpha) = w(\beta) = w(\gamma)$.
\end{definition}

\begin{definition}[$A_{\alpha,\beta,\gamma}$, $\chi_i^m$]
Let $(\alpha,\beta,\gamma)$ be a well-positioned triplet of nodes with $\alpha\in L^i_n$
and $w(\alpha) = U$. 
For each $m\in \N$ with $m \geq \lev(\alpha)$, define two nodes $\chi_0^m$ and $\chi_1^m$ as follows:

$\chi_0^{\lev(\alpha)} = \alpha$, $\chi_1^{\lev(\alpha)} = \beta$,  $\chi_1^{\lev(\alpha)+1} = \gamma$.

$$\chi_0^{m+1} = \node(\{\chi_0^m,\chi_1^m\})$$

$$\chi_1^{m+2} = \node(\{\chi_0^m,\chi_1^{m+1}\})$$

Note that $\chi^m_j\in \lev_{m}$. 


Let $A_{\alpha,\beta,\gamma} = \{\rho \mid \rho \geq \alpha \text{ or } \gamma\}
\cup \{\chi_i^m \mid m\in \N, i\in \{0,1\}\}$.
\end{definition}
The sets of the form $A_{\alpha,\beta,\gamma}$ will turn out to be exactly the
sets in $\{k(\phi) \mid \phi\in J_{2,n}$.
\begin{proposition} For each well-positioned triplet $(\alpha,\beta,\gamma)$ of nodes,
there is a formula
$\phi^{\alpha,\beta,\gamma}$ such that $k(\phi^{\alpha,\beta,\gamma}) = A_{\alpha,\beta,\gamma}$.
\end{proposition}
\begin{proof}
Let $A = A_{\alpha,\beta,\gamma}$. Let $S = \{\rho \mid \lev(\rho)\leq \lev(\alpha) + 2\text{ and } \rho \notin A\}$.

Let $\phi = \phi^{\alpha,\beta,\gamma}$ be 
$\psi_0\wedge \psi_1\wedge \psi_2 \wedge \psi_3$ where
$\psi_0$ is
\begin{displaymath}
\bigwedge_{\rho\in S}\phi'_\rho,
\end{displaymath}
$\psi_1$ is
\begin{displaymath}
\neg\neg(\phi_{\chi_0^{\lev(\alpha) + 2}}\vee\phi_{\chi_1^{\lev(\alpha)+2}}),
\end{displaymath}
$\psi_2$ is
\begin{displaymath}
\bigwedge_{x_i \in U} x_i,
\end{displaymath}
and $\psi_3$ is
\begin{displaymath}
\bigwedge_{x_i\notin U} (x_i \rightarrow (\phi_{\chi_0^{\lev(\alpha)+2}}\vee \phi_{\chi_1^{\lev(\alpha) + 2}}))
\end{displaymath}

We first show that $A\subset k(\phi)$.
Note that $A$ is upward closed, and that since $A$ is upward closed, we automatically have that
$\mu\Vdash \psi_0$ for all $\mu\in A$.

First, let $\mu\in A$ satisfy $\mu \geq \chi_0^{\lev(\alpha) + 2}$ or
$\mu\geq \chi_1^{\lev(\alpha)+2}$ . Then immediately we have $\mu\Vdash \psi_1$ and $\mu\Vdash \psi_3$.  Since
$\mu$ is above something which forces every atomic formula in $U$, we have that $\mu\Vdash \psi_2$.

Now, let $\mu\in A$ be $\chi_0^{\lev(\alpha) + m}$ or $\chi_1^{\lev(\alpha) + m}$ for some $m > 2$.  Since all of $\mu$'s
immediate successors force $\psi_1$ (by induction), and $\psi_1$ is a doubly-negated formula,
$\mu\Vdash \psi_1$. By construction, $\mu\Vdash \psi_2$.  Similarly, as the only nodes in $A$
which force atomic formulas not in $U$ must be above $\chi_0^{\lev(\alpha) + 2}$ or $\chi_1^{\lev(\alpha) + 2}$, $\mu\Vdash \psi_3$.

We now show that $k(\phi)\subset A$. Let $\mu\Vdash \phi$.  If
$\lev(\mu) \leq \lev(\alpha) + 2$ then the fact that $\mu\Vdash \psi_0$ implies
that $\mu\in A$.  

Let $\lev(\mu) > \lev(\alpha) + 2$ and assume (for a proof by induction) that we have proven that 
$k(\phi)\cap (\bigcup_{1\leq i \leq\lev(\mu)-1}L^i_n) = A\cap (\bigcup_{1\leq i \leq\lev(\mu)-1}L^i_n)$.  By
the construction of $K_n$, $\mu$ must be some $\node(T,U')$.  $T$
must contain some element of $A\cap L_n^{\lev(\mu) - 1}$, which, by induction,
must be either $\chi_0^{\lev(\mu) - 1}$ or $\chi_1^{\lev(\mu) - 1}$.  The choices for the set of minimal elements of $T$ are thus either
$\{\chi_0^{\lev(\mu) - 1},\chi_1^{\lev(\mu) - 1}\}$ or
$\{\chi_0^{\lev(\mu) - 2},\chi_1^{\lev(\mu) - 1}\}$.  The fact that $\mu\Vdash \psi_2$
forces $U'$ to equal $U$, and $\mu$ must be either $\chi_0^{\lev(\mu)}$ or
$\chi_1^{\lev(\mu)}$.
\end{proof}

The fact that these $\phi^{\alpha,\beta,\gamma}$ are in $J_{2,n}$ is seen by observing that
for sufficiently large $m$, $|A_{\alpha,\beta,\gamma}\cap L^m_n| = 2$,
that every $\rho\in A_{\alpha,\beta,\gamma}$ has some $\rho' < \rho$
such that $\rho'\in A_{\alpha,\beta,\gamma}$, 
and the following simple lemma.
\begin{lemma}\label{areminimal} Let $\phi$ be such that there is an $m$ such that
$|k(\phi)\cap L^m_n| = 1$.  Then for all $m' \geq m + 2$,
$|k(\phi)\cap L^m_n| = 0$.  In particular,
$|k(\phi)| < \infty$.
\end{lemma}
\straightforward

We will now observe that there are infinitely many distinct sets of the form $A_{\alpha,\beta,\gamma}$ in $K_n$
if $n\geq 2$.  Let $m$ be an arbitrary natural number, and let $|L^i_n| \geq 3m$. Find $m$ disjoint triples
$(\alpha_j,\beta_j,\delta_j)$ in $L^i_n$.  Then each $(\alpha_j,\beta_j,\node(\{\beta_j,\delta_j\}))$ is
well-positioned and the sets $A_{\alpha_j,\beta_j,\node(\{\beta_j,\delta_j\})}$ are pairwise unequal.  Thus
there are more than $m$ distinct sets of the form $A_{\alpha,\beta,\gamma}$ and $m$ was arbitrary.
 
We now show that there are no other minimal formulas and that every element of $J_n-J_{1,n}$ has a minimal element
of $J_n-J_{1,n}$ below it with the following proposition.

\begin{proposition} Let $\phi$ be a join-irreducible such that $k(\phi)$ is infinite. Then there is a
well-posititioned triplet
$(\alpha,\beta,\gamma)\in k(\phi)$ such that $A_{\alpha,\beta,\gamma}\subset k(\phi)$. 
\end{proposition}
\begin{proof}
By Lemma~\ref{hypfails}, if we find a well-positioned tuple $(\alpha,\beta,\gamma)$ such
that $\type_\phi(\alpha) = \type_\phi(\beta) = \type_\phi(\gamma)$ we will
be done.

Let $\alpha_0$ be such that $\type_\phi(\alpha_0) = \mintype(\phi)$. Since $k(\phi)$ is
infinite, there must be some $\beta_0\ne\alpha_0$ such that $\lev(\beta_0) = \lev(\alpha_0)$.

Let $\alpha_1 = \node(\{\beta_0,\alpha_0\})$. We have that $\type_\phi(\alpha_1) = \mintype(\phi)$.
Since $k(\phi)$ is infinite, there must be a $\beta_1\ne \alpha_1$ such that
$\lev(\beta_1) = \lev(\alpha_1)$

First assume that $\beta_1 < \alpha_0$. If there is any $\gamma'$ such that 
$\lev(\gamma') < \lev(\alpha_2)$ and $\gamma' \not \geq \alpha_2$ then we may
take $\alpha = \beta_2$, $\beta=\alpha_2$, $\gamma = \node(\{\alpha_2,\gamma'\})$. Otherwise,
there must be some $\gamma'$ such that $\lev(\gamma') < \lev(\alpha_2)$ and $\gamma' \not \geq \beta_2$ (otherwise $\alpha_2$ would equal $\beta_2$).  Then we may take $\alpha = \alpha_2$,
$\beta = \beta_2$ and $\gamma = \node(\{\beta_2,\gamma'\})$.

If $\beta_1\not <\alpha_0$ then repeat the argument of the above paragraph with $\beta_2 = \node(\{\alpha_0,\beta_1\})$ and $\alpha_2 = \node(\{\beta_1,\alpha_1\})$.
\end{proof}
\end{proof}
\begin{corollary}\label{threej3}
For any join-irreducible $\phi$, the following are equivalent:

1. $\phi\in J_{3,n}$.

2. There are three incomparable nodes $\alpha_1$, $\alpha_2$, $\alpha_3$ such 
that for $1\leq i \leq 3$, $\type_\phi(\alpha_i) = \mintype(\phi)$.

3. For any $r$, there is an $m$ such that $|\type_{\phi}^{-1}(\mintype(\phi))\cap L^m_n| \geq r$.
\end{corollary}
\begin{proof}
$1\implies 2$: Let $\phi\in J_{3,n}$
and let 
$(\alpha,\beta,\gamma)$ be a well-positioned tuple of nodes in $k(\phi)$.

By assumption, there is a $\delta\in k(\phi) - A_{\alpha,\beta,\gamma}$.
Let $\type_\phi(\mu)=\mintype(\phi)$. First assume that $\mu\in A_{\alpha,\beta,\gamma}$.

Then all $\nu\in A_{\alpha,\beta,\gamma}$ of sufficiently high level must satisfy
$\type_\phi(\nu) = \mintype(\phi)$.  Since $k(\phi)$ is infinite, every element of $k(\phi)$
must have a predecessor in $k(\phi)$: if some $\mu'$ didn't, then $\phi$ would not be
join-irreducible as it would be equivalent to $(\phi\wedge \phi'_{\mu'})\vee\phi_{\mu'}$.

Thus, we can find some predecessor $\delta'$ of $\delta$ on the same level as two elements
$\chi_0$, $\chi_1$
of $A_{\alpha,\beta,\gamma}$ satisfying $\type_\phi(\chi_0) = \type_\phi(\chi_1) = \mintype(\phi)$.
Then we can take our three incomparable nodes to be $\node(\{\delta',\chi_0\})$, $\node(\{\delta',\chi_1\})$, and $\node(\{\delta',\chi_0,\chi_1\})$.

Now assume that $\mu\notin A_{\alpha,\beta,\gamma}$.  Take a predecessor $\mu'$ of $\mu$ that
is on the same level as two elements $\chi_0$, $\chi_1$ of $A_{\alpha,\beta,\gamma}$. Then we may take our three incomparable nodes to be $\node(\{\mu',\chi_0\})$, $\node(\{\mu',\chi_1\})$ and
$\node(\{\mu',\chi_0,\chi_1\})$.

$2\implies 3$: Let $\alpha_1$, $\alpha_2$, and $\alpha_3$ be three incomparable nodes with
$\type_\phi(\alpha_i) = \mintype(\phi)$ for any $i$ between 1 and 3. Suppose the maximum level
of the $\alpha_i$'s is $m$.

Then at level $m+1$ we have four nodes of type $\mintype(\phi)$: $\node(\{\alpha_1,\alpha_2,\alpha_3\})$, $\node(\{\alpha_1,\alpha_2\})$, $\node(\{\alpha_2,\alpha_3\})$, $\node(\{\alpha_1,\alpha_3\})$.  Clearly, if at any level $m'$, we have $r$ nodes of type $\mintype(\phi)$, then at level $m'+1$, we have
at least $(r\text{ choose }2)$ nodes of type $\mintype(\phi)$.  Since the function 
$r\mapsto (r\text{ choose }2)$ is strictly increasing for $r \geq 4$, we are done.

$3 \implies 1$: Clearly $k(\phi)$ is infinite, so $\phi\notin J_{1,n}$. If $\phi\in J_{2,n}$,
then $|k(\phi)\cap L^m_n|$ would be 2 for all sufficiently large $m$.
\end{proof}
\section{$J_{3,n}$}

We begin with some facts about qsls.

\begin{lemma}\label{finite} Any finite set of elements in a bqsl has only finitely
many maximal lower bounds.
\end{lemma}
\begin{proof}
We will prove it for a set of three elements. The general case is by induction.

Let $p$, $q$, and $r$ be elements of a bqsl.  Let $\{s_1,\ldots,s_m\}$ be the
maximal lower bounds of $p$ and $q$.  For each $i$ from 1 to $m$, let $T_i$
be the set of maximal lower bounds of $r$ and $s_i$. 
Then $\bigcup_i T_i$ is finite and every element less than all three of $p$,
$q$ and $r$ is less than some element of $\bigcup_i T_i$: since it's less than
$p$ and $q$, it's less than some $s_i$ and therefore is less than some element of
$T_i$ since it's less than $r$.  The conclusion follows
\end{proof}
\begin{lemma}\label{oneelement} If $Q^*\subset Q^{**}$ are finite bqsl's, then there
is a sequence $Q^* = Q_0\subset Q_1\subset\cdots \subset Q_m = Q^{**}$ such that
each $Q_i$ is a bqsl, $|Q_{i+1}| = |Q_i| + 1$, and each inclusion of $Q_i$ in 
$Q_{i+1}$ is a bqsl embedding.
\end{lemma}
\begin{proof}
Let $q$ be a minimal element of $Q^{**} -Q^*$. Let $Q_1 = Q^*\cup \{q\}$.  This is a bqsl:
for any $p\in Q^*$, the maximal lower bounds of $p$ and $q$ must be in $Q^{**}$, therefore
they must be in $Q^*\cup \{q\}$ since $q$ was minimal in $Q^{**}-Q^*$. Repeat the process.
\end{proof}

\begin{lemma}\label{incomp} Let $\phi\in J_{3,n}$. Let $S$ be a finite collection of join-irreducible
formulas such that for all $\psi\in S$, $\phi\not\leq \psi$.  There is a $\rho\in J_{3,n}$
such that $\rho < \phi$ and for all $\psi\in S$, $k(\rho)\cap k(\psi)$ is finite.
\end{lemma}
\begin{proof}
Let $\alpha$ be such that $\type_\phi(\alpha)=\mintype(\phi)$. Let $S = \{\psi_1,\ldots,
\psi_q\}$. For $1\leq i \leq q$, let $\alpha_i\in k(\phi) - k(\psi_i)$
By Corollary~\ref{threej3}, find an $m$ greater than the level of each $\psi_i$ such 
that there are four nodes $\gamma_1,\ldots,\gamma_4$ in $L^m_n$ such that
for $1\leq i \leq 4$, $\type_\phi(\gamma_i) = \mintype(\phi)$.

Let $\delta_j = \node(\{\alpha_i\mid 1\leq i \leq q\}\cup\{\gamma_j\})$ for $1\leq j \leq 3$.
Let $m$ be the maximum level of a $\delta_i$. Let $T = \{\beta\mid \lev(\beta)\leq m,
(\forall i)\,\beta\not\geq\delta_i\}$.   Then a $\rho$ satisfying the required properties is
\[ \neg\neg(\phi_{\delta_1}\vee\phi_{\delta_2}\vee\phi_{\delta_3})\wedge \bigwedge_{\beta\in T}\phi'_{\beta}.\]

We will first show that $k(\rho)\subset k(\phi)$. This will be done by induction on the level
of nodes in $k(\rho)$.  First, let $\alpha\in k(\rho)$ be such that $\lev(\alpha)\leq m$. Then $\alpha$ must be in $k(\phi)$ since
it forces $\bigwedge_{\beta\in T}\phi'_{\beta}$ and is thus above some $\delta_i$, which forces
$\phi$.

Now consider a node of the form $\alpha =\node(S,U)$, where $\lev(\alpha) > m$. Since 
$k(\rho)\cap L^m_n = \{\delta_1,\delta_2,\delta_3\}$, $\alpha$ must be below some
$\delta_i$.  Since all elements of $S$ force $\phi$ by induction, and $\type_\phi(\delta_i) = \mintype(\phi)$, $\alpha$ forces $\phi$.

Note that $k(\rho)\ne k(\phi)$ as $\gamma_4\in k(\phi)-k(\rho)$.

Finally, since no $\delta_i$ forces any $\psi_j$, $k(\rho)\cap k(\psi_j)\cap L^m_n = \emptyset$.
Thus $k(\rho)\cap k(\psi_j)$ is finite.
\end{proof}

\begin{proposition}\label{j3} The set $J_{3,n}$ is isomorphic to
$Q'$.
\end{proposition}
\begin{proof}
Let $J$ be $J_{3,n}$ with a minimum element $\bot$ added.

First we observe that $J$ is a bounded locally finite quasisemilattice. It has a minimum element.  Given
any two formulas $\phi$ and $\psi$, $\phi\wedge\psi$ is their greatest lower bound
in $H_n$.  Although it may not be join-irreducible, it can be written as 
$\rho_1\vee\cdots\vee \rho_m$ with each $\rho_i$ join-irreducible.  The maximal elements
among those $\rho_i$ in $J_{3,n}$ are then the maximal lower bounds of $\phi$ and $\psi$ in $J_{3,n}$ and there are
only finitely many of them. If no $\rho_i$ is in $J_{3,n}$, then the maximum lower bound of $\phi$ and $\psi$ in $J$ is $\bot$.

We now show that $J$ is locally finite.  Given join-irreducible $\phi$ and $\psi$,
their maximal lower bounds are the maximal join-irreducibles less than $\phi\wedge\psi$.
By Corollary~\ref{subformulas} and Lemma~\ref{generaljoin} these are formed out of subformulas
of $\phi$ and $\psi$ by $\wedge$ and $\vee$. Iterating the process still yields formulas
formed out of subformulas of $\phi$ and $\psi$ by $\wedge$ and $\vee$.  Thus, there can only be finitely many such formulas.

We now prove that it is the universal countable homogeneous locally finite bounded quasisemilattice by
showing that for any finite bounded quasisemilattice $Q_1$, any bqsl embedding $f$ from $Q_1$ into
$J_n$, and any bqsl embedding $g$ from $Q_1$ into $Q_2$ where $|Q_2| = |Q_1| + 1$, there
is an extension $h$ of $f$ along $g$ from $Q_2$ into $J_n$ which is also a bqsl embedding.

Let $q$ be the unique element of $Q_2 - Q_1$. Let $U = \{f(p) \mid p\in Q_1, p > q\}$,
$K = \{f(p)\mid p\in Q_1, p \not< q,p\not>q \}$ and $L = \{f(p)\mid p\in Q_1, p < q\}$.

Assume $U$ is nonempty. Then by Lemma~\ref{finite}, there must be a minimum element of $U$.
Let it be $u^*$. If $U$ is empty, let $u^*$ be $\top$.

First assume that $L$ has more than one maximal element.  

Since $u^*$ is not less
than any element of $K$, we can find, for each $\psi\in K\cup L$, an element $\alpha_\psi\in k(u^*)$
that is not in $k(\psi)$.  
By Corollary~\ref{threej3}, we can find two incomparable nodes
$\beta_1$ and $\beta_2$ at a level greater than any $\psi\in K\cup L$. 
Let $\beta = \node(\{\alpha_\psi\mid\psi\in K\}\cup\{\beta_1\})$ and let
$\beta' = \node(\{\alpha_\psi\mid\psi\in K\}\cup\{\beta_2\})$.  Note that $\beta$
and $\beta'$ are incomparable,
$\type_{u^*}(\beta)=\type_{u^*}(\beta')=\mintype(u^*)$ and
for all $\psi\in K\cup L$, $\beta,\beta'\notin k(\psi)$.


Our new element $\phi$ (representing $q$) will be 
$\phi'_\beta \rightarrow (\bigvee_{\rho\in L} \rho\vee\phi_\beta)$.

We will denote $\bigvee_{\rho\in L} \rho$ by $\bigvee L$ from now on.

Let $R_1 = \{\beta\}$.  For $n > 1$, let
$R_n = \{\node(R)\mid R\subset \bigcup_{i < n} R_i\cup k(\bigvee L) \wedge
R\cap \bigcup_{i<n} R_i \ne \emptyset
\}$.

\begin{lemma}\label{decomposition} $k(\phi) = \bigcup_i R_i \cup k(\phi_\beta)\cup k(\bigvee L)$.
\end{lemma}
\begin{proof}
We will first show that $\bigcup_i R_i \cup k(\phi_\beta)\cup k(\bigvee L)\subset k(\phi)$.

Clearly $k(\phi_\beta)$, $k(\bigvee L)$, and $R_1$ are subsets of $k(\phi)$.

We will also show that no element of any $R_i$ forces $\phi'_\beta$.  Clearly the sole
element of $R_1$ does not force it.

Assume $R_n\subset k(\phi)$ and that no element of $R_i$ for
$i \leq n$ forces $\phi'_\beta$. 
Let $R\subset\bigcup_{i < n} R_i\cup k(\bigvee L)$
and $R\cap \bigcup_{i<n} R_i \ne \emptyset$.  

Since there is an element of $R$ that does not force $\phi'_\beta$, $\node(R)$ doesn't force it either,
and thus forces $\phi'_\beta \rightarrow (\bigvee_{\rho\in L} \rho\vee\phi_\beta))$
since all of its successors force $\phi$.  Thus, $\node(R)$ forces $\phi$.

We will now show that $k(\phi)\subset \bigcup_i R_i\cup k(\phi_\beta)\cup k(\bigvee L)$.
Suppose $\alpha \in k(\phi)$, $\alpha\notin k(\phi_\beta)$, and
$\alpha\notin k(\bigvee L)$.  

Since $\alpha\not\Vdash \phi_\beta\vee\bigvee L$, we must have $\alpha\not\Vdash \phi'_\beta$.
Thus $\alpha \leq \beta$. It follows that every node in $\bigcup_i R_i$ is less than $\beta$

Given any $\alpha'\leq \beta$, let $\alpha' = \alpha'_0 < \alpha'_1 < \cdots < \alpha'_{n(\alpha')} = \beta$ where $\alpha'_{i+1}$ is an immediate successor of $\alpha'_i$ and $n(\alpha') $ is as large as possible.

We will show by induction on $m$ that for all $\alpha'\leq \beta$ such that
$\alpha'\Vdash \phi$, $\alpha'\in R_{m+1}$
iff $n(\alpha') = m$.

The case where $m = 0$ is clear.

Suppose that it's true for $m$ and we'll show it true for $m+1$. First, let
$\alpha'\in R_{m+1}$.  Then $\alpha' = \node(R)$, where $R\subset \bigcup_{i\leq m} R_i\cup
k(\phi_\beta)\cup k(\bigvee L)$.  The maximum distance to $\beta$ is given by the maximum distance
from one if its immediate successors plus one.  

Conversely, if $n(\alpha')\leq m$, then since all of its successors must force $\phi$,
the ones less than $\beta$ are in $R_m$, and the ones not less than
$\beta$ must be in $k(\phi_\beta)\cup k(\bigvee L)$, $\alpha'$ must be in $R_{m+1}$.

Thus we are done, as each $\alpha' \leq \beta$ such that $\alpha \Vdash \phi$ is
in some $R_m$ (namely, $R_{n(\alpha')}$).
\end{proof}

We must show that $\phi$ is different from every element of $L$, $K$, and $U$, that it is less
than every element of $u^*$, greater than every element of $L$, incomparable with every element
of $K$, and for each $\psi\in K$, the maximal lower bounds of $\phi$ and $\psi$ are in $L$.

Clearly, $\phi$ is above every element of $L$.  Thus, $\phi$ is different from every element of
$L$, since it $L$ is assumed to have more than one maximal element.

We have that $\phi$ is below $u^*$ since $\phi_\beta$ and $\bigvee L$ are
and by induction on $i$ the type of each element of $R_i$ is $\mintype(u^*)$. We have that $\phi$ is different from $u^*$ as
$\beta'\Vdash \phi'_\beta$ but $\beta'\not\Vdash \bigvee L$, so $\beta'\in k(u^*)-k(\phi)$.  

We will show that $\phi$ is incomparable with each element of $K$. Fix a $\psi\in K$. Since 
$\beta\in k(\phi)$ but not in $k(\psi)$, we have $\psi\not\leq \phi$. On the other hand, the intersection of $k(\psi)$ with $k(\phi)$ must be contained in $k(\bigvee L)\cup k(\phi_\beta)$,
as every element of $k(\phi)$ is either below $\beta$, or contained in $k(\bigvee L)\cup k(\phi_\beta)$ by Lemma~\ref{decomposition}. Since
$\psi$ cannot be below $\bigvee L\vee\phi_\beta$ by Lemma~\ref{generaljoin},  there must be an element
$k(\psi)$ not in $k(\phi)$. Of course, it follows that $\phi$ is not equal to any element of $K$.

Let $\psi\in K$. In order to show that $\psi$ and $\phi$ have the same maximal lower bounds
in $K\cup L\cup U\cup \{\phi\}$ as they do in $J$, we must show that every $\chi\in J_{3,n}$ such that 
$\chi\leq \psi\wedge \phi$ is less than some element of $L$. But, as observed above
$k(\psi)\cap k(\phi)$ is contained in $k(\bigvee L\vee\phi_\beta)$. Since $\chi$ is join-irreducible and $k(\chi)$ is infinite,
$\chi$ is less than some $\rho\in L$.

This completes the case where $L$ has more than one maximal element. If $L$ has just one maximal element and it is not $\bot$, then apply Lemma~\ref{incomp} to get a $\chi$ such that: $\chi$ is below
$u^*$, no $\psi\in K$ is below $\chi$, and if $\chi'\in J$ is below $\chi$ and $\psi$, for
any $\psi\in K\cup L$ then $\chi'=\bot$. Add $\chi$ to $L$ and
proceed as before.

Finally, suppose $L = \{\bot\}$. Then apply Lemma~\ref{incomp} to get a $\chi$
below $u^*$ and such that $k(\chi\wedge\psi)$ is finite (and hence the greatest
lower bound in $J$ is $\bot$) for each $\psi\in K$.
\end{proof}

\end{document}